\documentclass[12pt]{article}

\begin{document}

\begin{center}
{\bf{GENERALIZED DIFFERENCE CES$\grave{A}$RO SEQUENCE SPACES OF FUZZY REAL
NUMBERS DEFINED BY ORLICZ FUNCTION}}
\end{center}

\vskip 0.5 cm

\centerline {\bf{B{\footnotesize inod} C{\footnotesize handra}
T{\footnotesize ripathy} $^*${\footnote {The work of the authors
was carried under University Grants Commision of India project
No.-{\it F}. No. $30-240/2004$ (RS) }} and S{\footnotesize tuti}
B{\footnotesize orgohain} $^{**}$}}

\vskip 0.5 cm

\noindent $^{*1}$Mathematical Sciences Division

\noindent Institute of Advanced Study in Science and
Technology

\noindent Paschim Boragaon; Garchuk; Guwahati-781035; Assam, {\bf{INDIA}}.

\noindent E-mail: $^*$tripathybc@yahoo.com;
tripathybc@rediffmail.com
\vskip 0.5 cm

\noindent $^{**}$ Department of Mathematics

\noindent Indian Institute of Technology, Bombay

\noindent Powai:400076, Mumbai, Maharashtra, {\bf{INDIA}}.

\noindent Email: $^{**}$stutiborgohain@yahoo.com
\vskip 2 cm

\noindent {\bf Abstract: }In this paper, we introduced different types of generalized difference Ces$\grave{a}$ro sequnecs spaces of fuzzy real
numbers defined by Orlicz function. We study some topological properties of these spaces. We obtain some inclusion relations involving these 
sequence spaces. These notions generalize many notions on difference Ces$\grave{a}$ro sequence spaces.\\



\vskip 0.5 cm

\noindent{\bf {1 Introduction}}

\vskip 1 cm

\noindent The concept of fuzzy set theory was introduced by L.A.Zadeh in 1965. Later on sequences of fuzzy numbers have been discussed by Esi (2006), Tripathy and Dutta (2007, 2010),
Tripathy and Sarma (2008a, 2008b) and many others.\\

\noindent {\bf{Definition 1.1}} Kizmaz (1981) defined the difference sequence spaces $\ell_\infty(\Delta), c(\Delta)$ and $c_0(\Delta)$ of crisp sets as follows:

$$Z(\Delta) = \{ (x_k) \in w : (\Delta x_k) \in Z \},$$

\noindent for $Z =\ell_\infty, c $ and $c_0$ and $\Delta x=(\Delta x_k) = (x_k - x_{k+1})$ for all $k \in N$.\\

The above spaces are Banach spaces, normed by,

$$ \Vert (x_k) \Vert_\Delta = \vert x_1 \vert + \displaystyle \sup_{k \geq 1} \vert \Delta x_k \vert.$$

The idea of Kizmaz (1981) was apllied to introduce different types of difference sequence spaces and study their different properties by Tripathy (2004),
Tripathy and Esi (2006), Tripathy and Mahanta (2004), Tripathy and Sarma (2008b), Tripathy, Altin and Et (2008) and many others.\\

Tripathy and Esi (2006) introduced the new type of difference sequence spaces, for fixed $m \in N$, as follows:

$$Z(\Delta_m) = \{ x=(x_k): (\Delta_m x_k) \in Z \},$$

\noindent for $Z =\ell_\infty, c $ and $c_0$ and $\Delta_m x=(\Delta_m x_k) = (x_k - x_{k+m})$ for all $k \in N$.\\

This generalizes the notion of difference sequence spaces studied by Kizmaz (1981).\\

The above spaces are Banach spaces, normed by,

$$\Vert (x_k) \Vert_{\Delta_m} = \displaystyle\sum_{r=1}^m\vert x_r \vert + \displaystyle \sup_{k \geq 1} \vert \Delta_m x_k \vert.$$

Tripathy, Esi and Tripathy (2005) further generalized this notion and introduced the following notion. For $m \geq 1$ and $n \geq 1$,

$$Z(\Delta_m^n) = \{x = (x_k): (\Delta_m^n x_k) \in Z \},\mbox{~for~} Z = \ell_\infty, c \mbox{~and~} c_0.$$

\noindent where $(\Delta_m^n x_k)=(\Delta_m^{n-1} x_k-\Delta_m^{n-1} x_{k+m})$, for all $k \in N$.\\

This generalized difference has the following binomial representation,

\begin{equation}
 \Delta_m^n x_k = \displaystyle \sum_{r=0}^n (-1)^r {n \choose r} x_{k+rm}.
\end{equation}

\noindent{\bf{Definition 1.2}} Ng and Lee (1978) defined the Ces$\grave{a}$ro sequence spaces $X_p$ of non-absolute type as follows:\\

$$x = (x_k) \in X_p \mbox{~if and only if~} \sigma(x) \in \ell_p, 1 \leq p < \infty,$$

\noindent where $\sigma(x)=\left(\frac{1}{n} \displaystyle\sum_{k=1}^n x_k \right)_{n=1}^\infty.$\\

Orhan (1983) defined the Ces$\grave{a}$ro difference sequence spaces $X_p(\Lambda)$, for $1 \leq p < \infty$ and studied their different 
properties and proved some inclusion results. He also obtained the duals of these sequence spaces.\\

Mursaleen, Gaur and Saifi (1996) defined the second difference Ces$\grave{a}$ro sequence spaces $X_p(\Delta^2)$, for $1 \leq p < \infty$ and studied 
their different topological properties and proved some inclusion results. They also calculated their dual sequence spaces.\\

Later on, Tripathy, Esi and Tripathy (2005) further introduced new types of difference Ces$\grave{a}$ro sequence spaces as $C_\infty(\Delta_m^n),
O_\infty(\Delta_m^n),C_p(\Delta_m^n),O_p(\Delta_m^n),$ and $\ell_p(\Delta_m^n)$, for $1 \leq p < \infty$.\\

For $m = 1$, the spaces $C_p(\Delta^n)$ and $C_\infty(\Delta^n)$ are studied by Et (1996-1997).\\ 
         
\noindent{\bf{Definition 1.3}} An Orlicz function is a function $M:[0,\infty) \rightarrow [0,\infty)$, which is continuous, non-decreasing and convex with 
$M(0) = 0, M(x) >0$, for $x>0$ and $M(x) \rightarrow \infty$, as $x \rightarrow \infty$.\\

An Orlicz function $M$ is said to satisfy $\Delta_2$-condition for all values of $x$, if there exists a constant $K>0, M(Lx) \leq KLM(x)$, for all 
$x>0$ and for $L>1$.\\

If the convexity of the Orlicz function is replaced by subadditivity i.e.$M(x + y) \leq M(x) + M(y)$, then this function is called as modulus function.\\

\noindent{\bf{Remark 1.1}} An Orlicz function satisfies the inequality $M(\lambda x) \leq \lambda M(x)$, for all $\lambda$ with $0 < \lambda < 1$.\\
        
Throughout the article $w^F,\ell^F,\ell_\infty^F$ represent the classes of {\it{all}}; {\it{absolutely summable}} and {\it{bounded}} sequences of fuzzy real numbers respectively.\\

\vskip 0.5 cm
 
\noindent{\bf{2 Definitions and background}}

\vskip 0.5 cm

\noindent{\bf{Definition 2.1}} A fuzzy real number $X$ is a fuzzy set on $R$ i.e. a mapping $X: R \rightarrow I (= [0,1])$ associating each real number $t$ with 
its grade of membership $X(t)$.\\

\noindent{\bf{Definition 2.2}} A fuzzy real number $X$ is called {\it{convex}} if $X(t) \geq X(s) \wedge X(r) = \mbox{min} (X(s), X(r))$, where $s < t < r$.\\

\noindent{\bf{Definition 2.3}} If there exists $t_0 \in R$ such that $X(t_0) =1$, then the fuzzy real number $X$ is called {\it{normal}}.\\
 
\noindent{\bf{Definition 2.4}} A fuzzy real number $X$ is said to be {\it{upper semi continuous}} if for each $\varepsilon>0,X^{-1}([0,a+\varepsilon))$, for all $a \in I$, is open in the usual topology of $R$.\\

The class of all {\it{upper semi-continuous, normal, convex}} fuzzy real numbers is denoted by $R(I)$.\\

\noindent{\bf{Definition 2.5}} For $X \in R(I)$, the $\alpha$-level set $X^\alpha$, for $0< \alpha \leq 1$ is defined by, $X^\alpha  = \{t \in R: X(t) \geq \alpha \}$. The $0$-level i.e. $X^0$ is the closure 
of strong $0$-cut, i.e. $X^0 = \mbox{cl}\{t \in R: X (t)>0\}$.\\

\noindent{\bf{Definition 2.6}} The absolute value of $X \in R(I)$ i.e. $\vert X \vert$ is defined by,

\[\vert X \vert (t)= \left\{
\begin{array}{l l}
\mbox{max} \{X(t),X(-t) \}, ~~~~ \mbox{for $t \geq 0$}\\
0\quad \mbox{otherwise}\\
\end{array} \right. \]

\noindent{\bf{Definition 2.7}} For $r \in R, \overline r \in R(I)$ is defined as,

\[\overline r(t) = \left\{
\begin{array}{l l}
1 ~~~~ \mbox{if $t=r$}\\
0\quad \mbox{if $t \neq r$}\\
\end{array} \right. \]

\noindent{\bf{Definition 2.8}} The additive identity and multiplicative identity of $R(I)$ are denoted by $\overline 0$ and $\overline 1$ respectively. The zero sequence of fuzzy real numbers is denoted by $\overline \theta$.\\

\noindent{\bf{Definition 2.9}} Let $D$ be the set of all closed bounded intervals $X = [X^L, X^R]$.\\ 

Define $d: D \times D \rightarrow R$ by $d(X,Y) = \mbox{max} \{ \vert X^L-Y^L \vert ,\vert X^R-Y^R \vert \}$. Then clearly $(D,d)$ is a complete 
metric space.\\

Define  by $\overline d:R(I) \times R(I) \rightarrow R(I)$ by $\overline d(X,Y)=\displaystyle\sup_{0 < \alpha \leq 1} d(X^\alpha,Y^\alpha)$, for $X,Y \in R(I)$. Then it is
well known that $(R(I),\overline d)$ is a complete metric space.\\

\noindent{\bf{Definition 2.10}} A sequence $X = (X_k)$ of fuzzy real numbers is said to converge to the fuzzy number $X_0$, if for every $\varepsilon >0$, there exists $k_0 \in N$ such that $\overline d(X_k,X_0) <\varepsilon$, for all $k \geq  k_0$.\\

\noindent{\bf{Definition 2.11}} A sequence space $E$ is said to be {\it{solid}} if $(Y_n) \in E$, whenever  $(X_n) \in E$ and $\vert Y_n \vert \leq \vert X_n \vert$,  for all $n \in N$.\\ 

\noindent{\bf{Definition 2.12}} Let $X = (X_n)$ be a sequence, then $S(X)$ denotes the set of all permutations of the elements of $(X_n)$ i.e. $S(X) = \{(X_{\pi(n)}):\pi \mbox{~is a permutation of {\it{N}}}\}$. A sequence space $E$ is said to be {\it{symmetric}} if $S(X) \subset  E$ for all $X \in E$.\\ 

\noindent{\bf{Definition 2.13}} A sequence space $E$ is said to be {\it{convergence-free}} if  $(Y_n) \in E$ whenever $(X_n) \in E$ and $X_n= \overline 0$ implies $Y_n= \overline 0$.\\

\noindent{\bf{Definition 2.14}} A sequence space $E$ is said to be monotone if $E$ contains the canonical pre-images of all its step spaces.\\ 

\noindent{\bf{Lemma 2.1}} A sequence space $E$ is solid implies that $E$ is monotone.\\ 

\noindent{\bf{Definition 2.15}} Lindenstrauss and Tzafriri (1971) used the notion of Orlicz function and introduced the sequence space:
                        
$$\ell_M = \left\{ x \in w: \displaystyle \sum_{k=1}^\infty M \left(\frac{\vert x_k \vert}{\rho}\right) , \mbox{~for some~} \rho>0\right\}$$

The space $\ell_M$ with the norm,

$$ \Vert x \Vert = \mbox{inf} \left\{ \rho > 0: \displaystyle \sum_{k=1}^\infty M \left(\frac{\vert x_k \vert}{\rho}\right) \leq 1 \right\}$$

becomes a Banach space, called an Orlicz sequence space. The space $\ell_M$ closely related to the space $\ell_p$ which is an Orlicz sequence space with $M(x) = x^p$, for $1 \leq p < \infty$.\\

Later on different classes of Orlicz sequence spaces were introduced and studied by Tripathy and Mahanta (2004), Et, Altin, Choudhary and Tripathy (2006), Tripathy, Altin and Et (2008), Tripathy and Sarma (2008a,2009,2011)many others.\\    

Let $m,n \geq 0$ be fixed integers and $1 \leq p < \infty$. In this article we introduced the following new types of generalized difference Ces$\grave{a}$ro sequence spaces of fuzzy real numbers:

$$C_p^F(M,\Delta_m^n)=\left\{(X_k):\displaystyle\sum_{i=1}^\infty\left(\frac{1}{i} \displaystyle\sum_{k=1}^i \left(M\left(\frac{\overline d(\Delta_m^n X_k, \overline 0)}{\rho}\right)\right)\right)^p<\infty,\mbox{~for some~} \rho>0\right\}.$$

$$C_\infty^F(M,\Delta_m^n)=\left\{(X_k):\displaystyle\sup_i \frac{1}{i}\left(\displaystyle\sum_{k=1}^i M\left(\frac{\overline d(\Delta_m^n X_k, \overline 0)}{\rho}\right)\right)<\infty,\mbox{~for some~} \rho>0\right\}.$$

$$\ell_p^F(M,\Delta_m^n)=\left\{(X_k):\displaystyle\sum_{k=1}^\infty\left(M\left(\frac{\overline d(\Delta_m^n X_k, \overline 0)}{\rho}\right)\right)^p<\infty,\mbox{~for some~} \rho>0\right\}.$$

$$O_p^F(M,\Delta_m^n)=\left\{(X_k):\displaystyle\sum_{i=1}^\infty \frac{1}{i} \left(\displaystyle\sum_{k=1}^i \left(M\left(\frac{\overline d(\Delta_m^n X_k, \overline 0)}{\rho}\right)\right)\right)^p<\infty,\mbox{~for some~} \rho>0\right\}.$$

$$O_\infty^F(M,\Delta_m^n)=\left\{(X_k):\displaystyle\sup_i \frac{1}{i} \displaystyle\sum_{k=1}^i M\left(\frac{\overline d(\Delta_m^n X_k, \overline 0)}{\rho}\right) <\infty,\mbox{~for some~} \rho>0\right\}.$$

\noindent{\bf{Lemma 2.2}} {\it{Let $1 \leq p < \infty$. Then}},\\

\noindent({\it{i}}) {\it{The space $C_p^F(M)$ is a complete metric space with the metric}},

$$\eta_1(X,Y)=\mbox{inf}\left\{\rho>0:\left(\displaystyle\sum_{i=1}^\infty \frac{1}{i} \displaystyle\sum_{k=1}^i \left(M\left(\frac{\overline d(X_k,Y_k)}{\rho}\right)\right)^p\right)^{\frac{1}{p}} \leq 1 \right\}.$$

\noindent({\it{ii}}) {\it{ The space $C_\infty^F(M)$ is a complete metric space with respect to the metric}},

$$\eta_2(X,Y)=\mbox{inf}\left\{\rho>0:\displaystyle\sup_i \frac{1}{i} \displaystyle\sum_{k=1}^i M\left(\frac{\overline d(X_k,Y_k)}{\rho}\right) \leq 1 \right\}.$$

\noindent({\it{iii}}) {\it{ The space $\ell_p^F(M)$ is a complete metric space with the metric}},

$$\eta_3(X,Y)=\mbox{inf}\left\{\rho>0:\left(\displaystyle\sum_{i=1}^\infty \left(M\left(\frac{\overline d(X_k,Y_k)}{\rho}\right)\right)^p\right)^{\frac{1}{p}} \leq 1 \right\}.$$

\noindent({\it{iv}}) {\it{The space $O_p^F(M)$ is a complete metric space with the metric}},

$$\eta_4(X,Y)=\mbox{inf}\left\{\rho>0:\left(\displaystyle\sum_{i=1}^\infty \frac{1}{i} \displaystyle\sum_{k=1}^i \left(M\left(\frac{\overline d(X_k,Y_k)}{\rho}\right)\right)^p\right)^{\frac{1}{p}} \leq 1 \right\}.$$

\noindent({\it{v}}) {\it{The space $O_\infty^F(M)$ is a complete metric space with respect to the metric}},

$$\eta_5(X,Y)=\mbox{inf}\left\{\rho>0:\displaystyle\sup_i \frac{1}{i} \displaystyle\sum_{k=1}^i M\left(\frac{\overline d(X_k,Y_k)}{\rho}\right) \leq 1 \right\}.$$

\noindent{\bf{Proof of lemma 2.2(i)}} Let $(X^{(u)})$ be a Cauchy sequence in $C_p^F(M)$ such that $X^{(u)} = (X_n^{(u)})_{n=1}^\infty$, for $i \in N$.\\

Let $\varepsilon >0$ be given. For a fixed $x_0>0$, choose $r>0$ such that $M\left(\frac{rx_0}{2}\right) \geq 1$. Then there exits a positive integer $n_0 = n_0 (\varepsilon)$ such that
	
$$\eta_1(X^{(u)},X^{(v)}) < \frac{\varepsilon}{rx_0}, \mbox{~for all~} u, v \geq n_0.$$

By the definition of $\eta_1$, we get:
	
\begin{equation}
 \mbox{inf}\left\{ \rho>0:\left(\displaystyle\sum_{i=1}^\infty \frac{1}{i}\displaystyle\sum_{k=1}^i\left(M\left(\frac{\overline d(X_k^{(u)},X_k^{(v)})}{\rho}\right)\right)^p\right)^{\frac{1}{p}} \leq 1 \right\} < \varepsilon, \mbox{~for all~} u,v \geq n_0.
\end{equation}

\noindent Which implies that,

\begin{equation}
 M\left(\frac{\overline d(X_k^{(u)},X_k^{(v)})}{\rho}\right) \leq 1, \mbox{~for all~} u,v \geq n_0.
\end{equation}

$\Rightarrow M\left(\frac{\overline d(X_k^{(u)},X_k^{(v)})}{\eta_1(X^{(i)},X^{(j)})}\right) \leq 1 \leq M\left(\frac{rx_0}{2}\right), \mbox{~for all~} u,v \geq n_0.$\\

Since $M$ is continuous, we get,\\

$\overline d(X_k^{(u)},X_k^{(v)}) \leq \frac{rx_0}{2}.\eta_1(X^{(u)},X^{(v)}),\mbox{~for all~} u,v \geq n_0.$\\

$\Rightarrow \overline d(X_k^{(u)},X_k^{(v)}) < \frac{rx_0}{2}. \frac{\varepsilon}{rx_0}=\frac{\varepsilon}{2}, \mbox{~for all~} u,v \geq n_0.$\\

$\Rightarrow \overline d(X_k^{(u)},X_k^{(v)}) < \frac{\varepsilon}{2}, \mbox{~for all~} u,v \geq n_0.$\\

Which implies that $(X_k^{(u)})$ is a Cauchy sequence in $R(I)$ and so it is convergent in $R(I)$ by the completeness property of $R(I)$.\\
         
Also, $\displaystyle\lim_u X_k^{(u)}=X_k$ , for each $k \in N$.\\

Now, taking  $v \rightarrow \infty$ and fixing $u$ and using the continuity of $M$, it follows from (3),

$$M\left(\frac{\overline d(X_k^{(u)},X_k)}{\rho}\right) \leq 1, \mbox{~for some~} \rho > 0.$$

Now on taking the infimum of such $\rho$’s, we get,

$$\mbox{inf}\left\{ \rho>0:\left(\displaystyle\sum_{i=1}^\infty \frac{1}{i}\displaystyle\sum_{k=1}^i\left(M\left(\frac{\overline d(X_k^{(u)},X_k)}{\rho}\right)\right)^p\right)^{\frac{1}{p}} \leq 1 \right\} < \varepsilon, \mbox{~for all~} u \geq n_0 \mbox{~(by (2))}.$$
         
Which implies that,
        
$$\eta_1(X^{(u)},X) < \varepsilon,\mbox{~for all~} u \geq n_0.$$
         
i.e. $\displaystyle\lim_u X^{(u)} =X$.\\

Now, we show that $X \in C_p^F(M)$.\\

We know that,\\
                           
$\overline d(X_k,\overline 0) \leq \overline d(X_k^{(u)},X_k) + \overline d(X_k^{(u)},\overline 0)$.\\

Since $M$ is continuous and non-decreasing, so we get,\\
         
\noindent $\displaystyle\sum_{i=1}^\infty \frac{1}{i} \displaystyle\sum_{k=1}^i \left(M\left(\frac{\overline d(X_k,\overline 0)}{\rho}\right)\right)^p$\\

\noindent $\leq \displaystyle\sum_{i=1}^\infty \frac{1}{i} \displaystyle\sum_{k=1}^i \left(M\left(\frac{\overline d(X_k^{(u)},X_k)}{\rho}\right)\right)^p +\displaystyle\sum_{i=1}^\infty \frac{1}{i} \displaystyle\sum_{k=1}^i \left(M\left(\frac{\overline d(X_k^{(u)},\overline 0)}{\rho}\right)\right)^p$\\ 

\noindent $< \infty$. (finite)\\

Which implies that $X \in C_p^F(M)$.\\

Hence $C_p^F(M)$ is a complete metric space.\\

This completes the proof.\\

\vskip 0.5 cm

\noindent{\bf{3 Main results}}

\vskip 0.5 cm
        
\noindent In this section, we prove the results relating to the introduced sequence spaces. The proof of the following result is a routine verification.\\

\noindent{\bf{Proposition 3.1}} {\it{The classes of sequences $C_\infty^F(M, \Delta_m^n),O_\infty^F(M, \Delta_m^n),C_p^F(M, \Delta_m^n),O_p^F(M, \Delta_m^n)$ and $\ell_p^F(M, \Delta_m^n)$, for $1 \leq p <\infty$, are metric spaces with respect to the metric,\\

\noindent $f(X,Y)= \displaystyle\sum_{k=1}^{mn} \overline d(X_k, \overline 0) + \eta(\Delta_m^n X_k, \Delta_m^n Y_k)$\\

\noindent where $Z = C_\infty^F,C_p^F,O_\infty^F,O_p^F,\ell_p^F$.}}

\vskip 0.5 cm

\noindent{\bf{Theorem 3.1}} {\it{Let $Z(M)$ be a complete metric space with respect to the metric $\eta$, the space $Z(M,\Delta_m^n)$ is a complete metric space with respect to the metric,

\noindent $f(X,Y)= \displaystyle\sum_{k=1}^{mn} \overline d(X_k, \overline 0) + \eta(\Delta_m^n X_k, \Delta_m^n Y_k)$\\

\noindent where $Z = C_\infty^F,C_p^F,O_\infty^F,O_p^F,\ell_p^F$.}}

\noindent{\it{Proof}} Let $(X^{(u)})$ be a Cauchy sequence in $Z(M,\Delta_m^n)$ such that $X^{(u)} =(X_n^{(u)})_{n=1}^\infty$.\\

We have for $\varepsilon >0$, there exists a positive integer $n_0 = n_0(\varepsilon)$ such that,

$$f(X^{(u)},X^{(v)}) < \varepsilon, \mbox{~for all~} u,v \geq n_0.$$

By the definition of $f$, we get:

\begin{equation}
\displaystyle\sum_{r=1}^{mn} \overline d(X_r^{(u)}, X_r^{(v)}) + \eta(\Delta_m^n X_k^{(u)},\Delta_m^n X_k^{(v)}) < \varepsilon,\mbox{~for all~} u,v \geq n_0.
\end{equation}
	
Which implies that, 

$$\displaystyle\sum_{r=1}^{mn} \overline d(X_r^{(u)}, X_r^{(v)}) < \varepsilon, \mbox{~for all~} u,v \geq n_0.$$

$$\Rightarrow \overline d(X_r^{(u)}, X_r^{(v)}) < \varepsilon, \mbox{~for all~} u,v \geq n_0, r=1,2,3...mn.$$

Hence $(X_r^{(u)})$ is a Cauchy sequence in $R(I)$, so it is convergent in $R(I)$, by the completeness property of $R(I)$, for $r =1,2,3 . . . . .mn$.\\

Let,
\begin{equation}
\displaystyle\lim_{u \rightarrow \infty} X_r^{(u)}=X_r,\mbox{~for~} r =1,2,3....mn
\end{equation}

Next we have,
\begin{equation}
\eta(\Delta_m^n X_k^{(u)},\Delta_m^n X_k^{(v)}) < \varepsilon,\mbox{~for all~} u,v \geq n_0
\end{equation}

Which implies that $(\Delta_m^n X_k^{(u)})$ is a Cauchy sequence in $Z(M)$, since $M$ is a continuous function and so it is convergent 
in $Z(M)$ by the completeness property of $Z(M)$.\\

Let, $\displaystyle\lim_u \Delta_m^n X_k^{(u)}=Y_k$ (say), in $Z(M)$, for each $k \in N$.\\
        
We have to prove that,

$$\displaystyle\lim_u X^{(u)}=X \mbox{~~and~~} X \in Z(M,\Delta_m^n).$$

For $k=1$, we have, from (1) and (5),
	
$$\displaystyle\lim_u X_{mn+1}^{(u)}=X_{mn+1},\mbox{~for~} m \geq 1, n \geq 1.$$

Proceeding in this way of induction, we get,

$$\displaystyle\lim_u X_k^{(u)}=X_k,\mbox{~for each~} k \in N.$$

Also, $\displaystyle\lim_u \Delta_m^n X_k^{(u)}= \Delta_m^n X_k,\mbox{~for each~} k \in N.$\\

Now, taking  $v \rightarrow \infty$ and fixing $u$  it follows from (4),

$$\displaystyle\sum_{r=1}^{mn} \overline d(X_r^{(u)}, X_r) + \eta(\Delta_m^n X_k^{(u)},\Delta_m^n X_k) < \varepsilon,\mbox{~for all~} u,v \geq n_0.$$
         
Which implies that,
        
$$ f(X^{(u)},X) <  \varepsilon,\mbox{~for all~} u \geq n_0.$$
         
i.e.$\displaystyle\lim_u X^{(u)}=X.$\\

Now, it is to show that $X \in Z(M,\Delta_m^n).$\\

We know that, 

\begin{eqnarray*}
 f(\Delta_m^n X_k, \overline 0) &\leq& f(\Delta_m^n X_k^{(i)}, \Delta_m^n X_k)+f(\Delta_m^n X_k^{(i)},\overline 0)\\
&<& \infty.
\end{eqnarray*}

Which implies that $X \in Z(M,\Delta_m^n).$\\

Hence $Z(M,\Delta_m^n)$ is a complete metric space.\\

This completes the proof of the theorem.\\
         
The proof of the following results is a consequence of the above result and lemma.\\
         
\noindent{\bf{Proposition 3.2}} {\it{Let $1 \leq p < \infty$.Then}},

\noindent ({\it{i}}) {\it{The space $C_p^F(M,\Delta_m^n)$ is a complete metric space with the metric}},

$$f_1(X,Y)=\displaystyle\sum_{r=1}^{mn} \overline d(X_r,Y_r) + \mbox{~inf}\left\{\rho>0: \left(\displaystyle\sum_{i=1}^\infty \frac{1}{i} \displaystyle\sum_{k=1}^i \left(M\left(\frac{\overline d(\Delta_m^n X_k, \Delta_m^n Y_k)}{\rho}\right)\right)^p\right)^\frac{1}{p} \leq 1 \right\}.$$
         
\noindent({\it{ii}}) {\it {The space $C_\infty^F(M,\Delta_m^n)$ is a complete metric space with respect to the metric}},

$$f_2(X,Y)=\displaystyle\sum_{r=1}^{mn} \overline d(X_r,Y_r) + \mbox{~inf}\left\{\rho>0: \displaystyle\sup_i \frac{1}{i}M\left(\frac{\overline d(\Delta_m^n X_k, \Delta_m^n Y_k)}{\rho}\right) \leq 1 \right\}.$$

\noindent({\it{iii}}) {\it{The space $\ell_p^F(M,\Delta_m^n)$ is a complete metric space with respect to the metric}},

$$f_3(X,Y)=\displaystyle\sum_{r=1}^{mn} \overline d(X_r,Y_r) +\mbox{inf}\left\{\rho>0:\left(\displaystyle\sum_{i=1}^\infty \left(M\left(\frac{\overline d(\Delta_m^n X_k,\Delta_m^n Y_k)}{\rho}\right)\right)^p\right)^{\frac{1}{p}} \leq 1 \right\}.$$

\noindent({\it{iv}}) {\it {The space $O_p^F(M,\Delta_m^n)$ is a complete metric space with the metric}},

$$f_4(X,Y)=\displaystyle\sum_{r=1}^{mn} \overline d(X_r,Y_r) + \mbox{inf}\left\{\rho>0:\left(\displaystyle\sum_{i=1}^\infty \frac{1}{i} \displaystyle\sum_{k=1}^i \left(M\left(\frac{\overline d(\Delta_m^n X_k, \Delta_m^n Y_k)}{\rho}\right)\right)^p\right)^{\frac{1}{p}} \leq 1 \right\}.$$
         
\noindent({\it{v}}) {\it{$O_\infty^F(M,\Delta_m^n)$ is a complete metric space with respect to the metric}},

$$f_5(X,Y)=\displaystyle\sum_{r=1}^{mn} \overline d(X_r,Y_r)+ \mbox{inf}\left\{\rho>0:\displaystyle\sup_i \frac{1}{i} \displaystyle\sum_{k=1}^i M\left(\frac{\overline d(\Delta_m^n X_k,\Delta_m^n Y_k)}{\rho}\right) \leq 1 \right\}.$$

\noindent{\bf{Theorem 3.2}} {\it{The classes of spaces $Z(M, \Delta_m^n)$ , where $Z = C_\infty^F, O_\infty^F, C_p^F ,O_p^F$ and $\ell_p^F$, for $1 \leq p < \infty$, are not monotone and as such are not solid for $m,n \geq 1$.}}\\

\noindent{\it{Proof.}}  Let us consider the proof for $C_p^F(M, \Delta_m^n)$. The proof follows from the following example:\\

\noindent{\it{Example 3.1}} Let $X_k = \overline k$, for all $k \in N$.\\

Let $m =3$ and $n =2$. Let $M(x) = \vert x \vert$, for all $x \in [0,\infty)$.\\ 

Then, we have, $\overline d(\Delta_3^2 X_k, \overline 0) = 0$, for all $k \in N$.\\

Hence, we get, for $1 \leq p< \infty$,

$$\left(\displaystyle\sum_{i=1}^\infty \frac{1}{i} \displaystyle\sum_{k=1}^i \left(M\left(\frac{\overline d(\Delta_3^2 X_k, \overline 0)}{\rho}\right)\right)^p\right)^\frac{1}{p} <\infty, \mbox{~for some~} \rho>0$$

Which implies that, $(X_k) \in C_p^F(M,\Delta_3^2)$.\\

Let $J = \{k : k \mbox{~is even~} \} \subseteq N$. Let $(Y_k)$  be the canonical pre-image of $(X_k)_J$ for the subsequence $J$ of $N$. Then,    

$$Y_k=
\cases{ X_k,~~ \mbox{for}~~   k \mbox{~odd} , \cr
        \overline 0,  ~~ \mbox{for}~~   k \mbox{~even}.} $$
     
But $(Y_k) \notin C_p^F(M,\Delta_3^2)$.\\
 
Hence the spaces are not monotone as such not solid.\\

This completes the proof.\\

\noindent {\bf{Remark 3.1}} For $m = 0$ or $n = 0$, the spaces $C_p^F(M)$ and $C_\infty^F(M)$ are neither solid nor monotone, where as the spaces $\ell_p^F(M), O_p^F(M)$ and $O_\infty^F(M)$ are solid and hence are monotone.\\

\noindent{\bf{Theorem 3.3}} {\it{The classes of spaces $Z(M, \Delta_m^n)$, where $Z = C_\infty^F, O_\infty^F,C_p^F, O_p^F$ and $\ell_p^F$, for $1 \leq  p < \infty$, are not symmetric, for $m, n \geq 1$}}.\\

\noindent{\it{Proof}}  Let us consider the proof for $C_\infty^F(M,\Delta_m^n)$. The proof follows from the following example:\\

\noindent{\it{Example 3.2}} Let $X_k = 1$, for all $k \in N$.\\

Let $m = 4$ and $n =1$. Let $M(x) = \vert x \vert$, for all $x \in [0,\infty)$.\\ 

Then, we have, $\overline d(\Delta_4 X_k, \overline 0) = 0$, for all $k \in N$.\\

Hence, we get,

$$\displaystyle\sup_i \frac{1}{i} \left(\displaystyle\sum_{k=1}^i M \left(\frac{\overline d(\Delta_4 X_k, \overline 0)}{\rho}\right)\right) < \infty, \mbox{~for some~} \rho>0.$$

Which implies that, $(X_k) \in C_\infty^F(M,\Delta_4)$.\\

Consider the rearranged sequence $(Y_k)$ of $(X_k)$ such that $(Y_k) = (X_1, X_2, X_4, X_3, X_9, X_5, X_{16}, X_6, X_{25},...)$ such that $\overline d(\Delta_4 Y_k, \overline 0) \approx k-(k-1)^2 \approx k^2$, for all $k \in N$.\\

Which shows,$$\displaystyle\sup_i \frac{1}{i} \left(\displaystyle\sum_{k=1}^i M \left(\frac{\overline d(\Delta_4 Y_k, \overline 0)}{\rho}\right)\right) = \infty, \mbox{~for some fixed~} \rho>0.$$

Hence, $(Y_k) \notin C_\infty
^F(M, \Delta_4)$ . It follows that the spaces are not symmetric.\\

This completes the proof.\\

\noindent{\bf{Theorem 3.4}}  {\it{The classes of spaces $Z(M, \Delta_m^n)$, where $Z = C_\infty^F, O_\infty^F,C_p^F, O_p^F$ and $\ell_p^F$, for $1 \leq  p < \infty$, are not convergence-free, for $m, n \geq 1$}}.\\

\noindent{\it{Proof}} Let us consider the proof for $C_\infty^F(M,\Delta_m^n)$. The proof follows from the following example:\\
         
\noindent{\it{Example 3.3}} Let $m = 3$ and $n =1$. Let $M(x) = x^3$, for all $x \in [0,\infty)$.\\

Consider the sequence $(X_k)$ defined as follows:\\

$$X_k(t)=
\cases{ 1+k^2 t,~~ \mbox{for}~~   t \in [-\frac{1}{k^2},0] , \cr
        1-k^2 t,~~ \mbox{for}~~   t \in [0, \frac{1}{k^2}] , \cr
        0, ~~~~~~~~~\mbox{otherwise}} $$
	         
Then,

$$ \Delta_3 X_k(t)=
\cases{ 1+\frac{k^2(k+3)^2}{2 k^2+6k+9} t,~~ \mbox{for}~~   t \in \left[-\frac{2 k^2+6k+9}{k^2(k+3)^2},0\right] , \cr
        1-\frac{k^2(k+3)^2}{2 k^2+6k+9}t,~~ \mbox{for}~~   t \in \left[0, \frac{2 k^2+6k+9}{k^2(k+3)^2}\right] , \cr
        0, ~~~~~~~~~\mbox{otherwise}} $$
          
Such that,$\overline d(\Delta_3 X_k, \overline 0)=\frac{2 k^2+6k+9}{k^2(k+3)^2}=\frac{1}{k^2}+\frac{1}{(k+3)^2}$.\\
         
We have,$\displaystyle\sup_i \frac{1}{i} \left(\displaystyle\sum_{k=1}^i M \left(\frac{\overline d(\Delta_3 X_k, \overline 0)}{\rho}\right)\right) < \infty, \mbox{~for some fixed~} \rho>0.$\\

Thus $(X_k) \in C_\infty^F(M,\Delta_3)$.\\

Now, let us take another sequence  $(Y_k)$ such that,

$$Y_k(t)=
\cases{ 1+\frac{t}{k^2},~~ \mbox{for}~~   t \in [-{k^2},0] , \cr
        1-\frac{t}{k^2},~~ \mbox{for}~~   t \in [0, {k^2}] , \cr
        0, ~~~~~~~~~\mbox{otherwise}} $$
                   
So that, 

$$ \Delta_3 Y_k(t)=
\cases{ 1+\frac{t}{2 k^2+6k+9},~~ \mbox{for}~~   t \in [-(2 k^2+6k+9),0] , \cr
        1-\frac{t}{2 k^2+6k+9},~~\mbox{for}~~   t \in [0, (2 k^2+6k+9)] , \cr
        0, ~~~~~~~~~\mbox{otherwise}} $$

But, $\overline d(\Delta_3 Y_k, \overline 0)=(2 k^2+6k+9)$ , for all $k \in N$.\\

Which implies that,$\displaystyle\sup_i \frac{1}{i} \left(\displaystyle\sum_{k=1}^i M \left(\frac{\overline d(\Delta_3 X_k, \overline 0)}{\rho}\right)\right) = \infty, \mbox{~for some fixed~} \rho>0.$\\

Thus, $(Y_k) \notin C_\infty^F(M,\Delta_3)$.\\

Hence $C_\infty^F(M,\Delta_m^n)$ is not convergence-free, in general.\\

This completes the proof.\\

\noindent{\bf{Theorem 3.5}}\\

\noindent ({\it{a}}){\it{$\ell_p^F(M,\Delta_m^n) \subset O_p^F(M,\Delta_m^n) \subset C_p^F(M,\Delta_m^n)$  and the inclusions are strict}}.\\

\noindent({\it{b}}){\it{$Z(M,\Delta_m^{n-1}) \subset Z(M,\Delta_m^n)$ (in general $Z(M,\Delta_m^i) \subset Z(M,\Delta_m^n)$, for $i = 1, 2, 3,...{n-1})$, for $Z =C_\infty^F,O_\infty^F,C_p^F ,O_p^F$   and $\ell_p^F$, for $1 \leq  p < \infty$}}.\\
       
\noindent({\it{c}}){\it{ $O_\infty^F(M,\Delta_m^n) \subset C_\infty^F(M,\Delta_m^n)$ and the inclusion is strict}}.\\

\noindent{\it{Proof}} ({\it{b}}) Let $(X_k) \in C_\infty^F(M,\Delta_m^{n-1})$ . Then we have, 

$$\displaystyle\sup_i \frac{1}{i} \left(\displaystyle\sum_{k=1}^i M \left(\frac{\overline d(\Delta_m^{n-1} X_k, \overline 0)}{\rho}\right)\right) < \infty, \mbox{~for some~} \rho>0.$$

Now we have, 
          
$\displaystyle\sup_i \frac{1}{i} \left(\displaystyle\sum_{k=1}^i M \left(\frac{\overline d(\Delta_m^n X_k, \overline 0)}{2\rho}\right)\right) = \displaystyle\sup_i \frac{1}{i} \left(\displaystyle\sum_{k=1}^i M \left(\frac{\overline d(\Delta_m^{n-1} X_k-\Delta_m^{n-1} X_{k+1}, \overline 0)}{2\rho}\right)\right) $\\

$\leq \displaystyle\sup_i \frac{1}{2}\left(\frac{1}{i} \left(\displaystyle\sum_{k=1}^i M \left(\frac{\overline d(\Delta_m^{n-1} X_k, \overline 0)}{\rho}\right)\right)\right)+
\displaystyle\sup_i \frac{1}{2}\left(\frac{1}{i} \left(\displaystyle\sum_{k=1}^i M \left(\frac{\overline d(\Delta_m^{n-1} X_{k+1}, \overline 0)}{\rho}\right)\right)\right) .$\\

$<\infty$.\\

Proceeding in this way, we have,
 
$Z(M,\Delta_m^i) \subset  Z(M,\Delta_m^n)$, for $0 \leq i < n$, for $Z = C_\infty^F, O_\infty^F,C_p^F, O_p^F$ and $\ell_p^F$, for $1 \leq  p < \infty$.\\

This completes the proof. \\
        
\noindent{\bf{Theorem 3.6}} ({\it{a}}) {\it{If $1 \leq p < q$, then}},\\

\noindent ({\it{i}}) $C_p^F(M,\Delta_m^n) \subset C_q^F(M,\Delta_m^n)$\\ 

\noindent ({\it{ii}}) $\ell_p^F(M,\Delta_m^n) \subset \ell_q^F(M,\Delta_m^n)$\\
 
\noindent ({\it{b}}) $C_p^F(M) \subset C_p^F(M,\Delta_m^n)$ , {\it{for all $m \geq 1$ and $n \geq 1$.}}\\

\noindent{\bf{Theorem 3.7}} {\it{Let $M, M_1$ and $M_2$ be Orlicz functions satisfying $\Delta_2$- condition. Then, for $Z = C_\infty^F, O_\infty^F,C_p^F, O_p^F$ and $\ell_p^F$, for $1 \leq  p < \infty$}},\\

\noindent({\it{i}})  $Z(M_1,\Delta_m^n) \subseteq Z(M \circ M_1,\Delta_m^n)$\\

\noindent({\it{ii}}) $Z(M_1,\Delta_m^n) \cap Z(M_2,\Delta_m^n) \subseteq Z(M_1+M_2,\Delta_m^n)$.\\
        
\noindent{\it{ Proof}} ({\it{i}}) Let $(X_k) \in Z(M_1,\Delta_m^n)$. For  $\varepsilon >0$, there exists $\eta>0$ such that $\varepsilon = M (\eta)$.\\

Then,$M_1\left(\frac{\overline d(\Delta_m^n X_k,L)}{\rho}\right)<\eta$, for some $\rho>0, L \in R(I)$.\\

Let $Y_k = M_1\left(\frac{\overline d(\Delta_m^n X_k, L)}{\rho}\right)$, for some $\rho>0, L \in R(I)$.\\

Since $M$ is continuous and non-decreasing, we get,

$$M(Y_k)=M\left(M_1\left(\frac{\overline d(\Delta_m^n X_k, L)}{\rho} \right)\right) <M(\eta)=\varepsilon,\mbox{~for some~} \rho>0.$$

Which implies that,  $(X_k) \in Z(M \circ M_1,\Delta_m^n)$.\\

This completes the proof.\\

({\it{ii}}) Let  $(X_k) \in Z (M_1,\Delta_m^n) \cap Z(M_2,\Delta_m^n)$.\\

Then,$M_1\left(\frac{\overline d(\Delta_m^n X_k, L)}{\rho}\right) < \varepsilon$ , for some $ \rho>0, L \in R(I)$.\\

and $M_2\left(\frac{\overline d(\Delta_m^n X_k, L)}{\rho}\right) < \varepsilon$ , for some $ \rho>0, L \in R(I)$.\\

The rest of the proof follows from the equality,

\begin{eqnarray*}
 (M_1+M_2)\left(\frac{\overline d(\Delta_m^n X_k, L)}{\rho}\right)&=& M_1\left(\frac{\overline d(\Delta_m^n X_k, L)}{\rho}\right)+M_2\left(\frac{\overline d(\Delta_m^n X_k, L)}{\rho}\right)\\
&<& \varepsilon+\varepsilon=2 \varepsilon, \mbox{~for some~} \rho>0.
\end{eqnarray*}
 
Which implies that $(X_k) \in Z(M_1+M_2,\Delta_m^n)$.\\	

This completes the proof.\\

\vskip 0.5 cm

\noindent{\bf{References}}

\vskip 0.5 cm

\noindent Altin Y, Et M, Tripathy BC (2004) The sequence space $\vert \overline N_p \vert (M,r,q,s)$ on seminormed
spaces. Applied Math \&\ Computation 154: 423-430.\\

\noindent Esi A (2006) On some new paranormed sequence spaces of fuzzy numbers defined by
Orlicz functions and statistical convergence. Mathematical Modelling and Analysis
11(4):379-388.\\

\noindent Et M (1996-1997) On some generalized Ces$\grave{a}$ro difference sequence spaces. Istanbul Univ
fen fak Mat Dergisi 55-56: 221-229.\\

\noindent Et M, Altin Y, Choudhary B, Tripathy BC (2006) Some classes of sequences defined by
sequences of Orlicz functions. Math Ineq Appl 9(2):335- 342.\\

\noindent Kizmaz H (1981) On certain sequence spaces. Canad Math Bull 24(2) :169-176.\\

\noindent Lindenstrauss J, Tzafriri L (1971) On Orlicz sequence spaces. Israel J Math 10: 379-390.\\

\noindent Mursaleen M, Gaur, Saifi AH (1996) Some new sequence spaces their duals and matrix
transformations. Bull Cal Math Soc 88 :207-212.\\

\noindent Ng PN, Lee PY (1978) Ces$\grave{a}$ro sequence spaces of non absolute type. Comment Math 20:
429-433.\\

\noindent Orhan C (1983) Ces$\grave{a}$ro difference sequence spaces and related matrix transformations.
Comm Fac Univ Ankara ser A 32 :55-63.\\

\noindent Shiue JS (1970) On the Ces$\grave{a}$ro sequence spaces. Tamkang Journal Math 1:19-25.\\
 
\noindent Tripathy BC (2004) On generalized difference paranormed statistically convergent
sequences. Indian J Pure Appl Math 35(5) :655-663.\\

\noindent Tripathy BC, Altin Y, Et M (2008) Generalized difference sequences spaces on
seminormed spaces defined by Orlicz functions. Math. Slovaca 58(3): 315-324.\\

\noindent Tripathy BC, Baruah A (2009), New type of difference sequence spaces of fuzzy real
numbers. Mathematical Modelling and Analysis, 14(3) : 391-397.\\

\noindent Tripathy BC, Borgohain S (2008) The sequence space $m(M,\phi,\Delta_m^n, p)^F$. Mathematical
Modelling and Analysis 13(4): 577-586.\\

\noindent Tripathy BC, Borgogain S (2011), Some classes of difference sequence spaces of fuzzy real
numbers defined by Orlicz function, Advances in Fuzzy Systems Article ID216414, 6 pages.\\

\noindent Tripathy BC, Dutta AJ (2007) On fuzzy real-valued double sequence spaces ....Mathematical and Computer Modelling 46(9-10):1294-1299.\\

\noindent Tripathy BC, Dutta AJ (2010), Bounded variation double sequence space of fuzzy real
numbers. Computers and Mathematics with Applications 59(2) : 1031-1037.\\

\noindent Tripathy BC, Esi A (2006) A new type of difference sequence spaces. Inter Jour Sci Tech
1(1) :11-14.\\

\noindent Tripathy BC, Esi A, Tripathy BK (2005) On a new type of generalized difference Ces$\grave{a}$ro
sequence spaces. Soochow J Math 31: 333-340.\\

\noindent Tripathy BC, Hazarika B (2011), $I$-convergent sequence spaces defined by Orlicz
functions. Acta Mathematica Applicata Sinica 27(1) : 149-154.\\

\noindent Tripathy BC, Mahanta S (2004) On a class of generalized lacunary difference sequence
spaces defined by Orlicz function. Acta Math Applicata Sinica 20(2):231-238.\\

\noindent Tripathy BC, Sarma B (2008a) Sequence spaces of fuzzy real numbers defined by Orlicz
functions. Math Slovaca 58(5): 621-628.\\

\noindent Tripathy BC, Sarma B, (2008b)Statistically convergent difference double sequence spaces.
Acta Mathematica Sinica 24(5): 737-742.\\

\noindent Tripathy BC, Sarma B (2009) Vector valued double sequence spaces defined by Orlicz
function. Mathematica Slovaca 59(6) : 767-776.\\

\noindent Tripathy BC, Sarma B (2011), Double sequence spaces of fuzzy numbers defined by Orlicz
functions. Acta Mathematicae Scientia 31B(1) : 134-140.\\

\noindent Tripathy BC, Sen M, Nath S (2011), $I$-convergence in probabilistic $n$-normed space. Soft
Computing DOI 10.1007/s0050-011-0799-8 (in press).

\end{document}